\documentclass[12pt]{article}

\usepackage{amsmath,amsthm}
\usepackage{latexsym}
\usepackage {pstricks,pstcol,pst-node}
\usepackage{color}
\usepackage{xspace}
\usepackage[colorlinks=true,
linkcolor=green,
filecolor=brown,
citecolor=green]{hyperref}
\def\red{\textcolor{red} }

\begin{document}

\begin{center}
{\Large
A Combinatorial Interpretation for a Super-Catalan Recurrence                           \\ 
}
\vspace{10mm}
DAVID CALLAN  \\
Department of Statistics  \\
University of Wisconsin-Madison  \\
1210 W. Dayton St   \\
Madison, WI \ 53706-1693  \\
{\bf callan@stat.wisc.edu}  \\
\vspace{5mm}
August 9, 2004
\end{center}

\vspace{5mm}

\centerline{   \textbf{Abstract}}
Nicholas Pippenger and Kristin Schleich have recently given a combinatorial interpretation for the
second-order super-Catalan numbers $(u_{n})_{n\ge 
0}=(3,2,3,6,14,36,...)$: they count  
``aligned cubic trees'' on $n$ internal vertices.  Here we give a combinatorial interpretation of
the recurrence $u_{n} =  
\sum_{k=0}^{n/2-1}\binom{n-2}{2k}2^{n-2-2k}u_{k}\,:$ it counts
these trees by number of deep interior vertices where deep interior
means ``neither a leaf nor adjacent to a leaf''.
\vspace{8mm}

{\Large \textbf{1 \quad Introduction}  }

For each integer $m\ge 1$, the numbers
\[
\frac{\binom{2m}{m} \binom{2n}{n} }{2\binom{m+n}{m}}
\]
satisfy the recurrence relation
\begin{equation}
u_{n}=\sum_{k\ge 0}2^{n-m-2k} \binom{n-m}{2k}u_{k}
\label{eq1}
\end{equation}
and hence are integers except when $m=n=0$ \cite{superballot}. For fixed $m$, we'll call them super-Catalan 
numbers of order $m$ (although other numbers go by this name too). 
For $m=0$ and $n\ge 1$, they are the odd central binomial 
coefficients $(1,3,10,35,\ldots)$ which count lattice paths of $n$ 
upsteps and $n-1$ downsteps. For $m=1$ and $n\ge 0$, they are the 
familiar Catalan numbers $(1,1,2,5,14,42,\ldots)$ with numerous 
combinatorial interpretations \cite[Ex 6.19]{ec2}. For $m=2$, three 
combinatorial interpretations have recently been given, in terms of 
(i) pairs of Dyck paths \cite{gesselxin04}, (ii) 
``blossom trees'' \cite{schaeffer03}, and (iii) ``aligned cubic trees'' 
\cite{pippenger}. The object of this note is to 
establish a combinatorial interpretation of the recurrence (\ref{eq1}) 
for $m=2$: 
it counts the just mentioned aligned cubic trees by number of vertices 
that are neither a leaf nor adjacent to a leaf.

For $m=1$, (\ref{eq1}) is known as Touchard's identity, and
in Section 2 we recall combinatorial interpretations of the recurrence 
for the cases $m=0,1$. In Section 3 we define aligned cubic trees and 
establish notation. In Section 4 we introduce configurations 
counted by the right side of (\ref{eq1}) for $m=2$, and in Section 5 we exhibit a 
bijection from them to size-$n$ aligned cubic trees.

\vspace*{3mm}

{\Large \textbf{2 \quad Recurrence for \emph{m}\,=\:0,\,1 }  }

For $m=0$, (\ref{eq1}) 
counts lattice paths of $n$ upsteps ($U$) and 
$n-1$ downsteps ($D$) by number $k$ of $DUU$s where, for example, 
$DDUUUDDUU$ has 2 $DUU$s. It also
counts  paths of $n$ $U$s and 
$n$ $D$s that start up by number $2k$ (necessarily even) of 
inclines ($UU$ or $DD$) at \emph{odd} 
locations. For example, $U_{\textrm{{\,\tiny 1\,}}}U_{\textrm{{\,\tiny 2\,}}}U_{\textrm{{\,\tiny 3\,}}}
D_{\textrm{{\,\tiny 4\,}}}D_{\textrm{{\,\tiny 5\,}}}D$ has four inclines, at 
locations $1,2,4$ and 5, but only the first and last are at odd locations. 

For $m=1$, (\ref{eq1}) counts Dyck $n$-paths (paths of $n\ U$s and $n\ 
D$s that
never dip below ground level) by number $k$ of $DUU$s. It also counts them 
by number $2k$ of inclines at \emph{even}
locations. See \cite{twobij04}, for example, for relevant bijections. 
Recall the standard ``walk-around'' bijection from full binary trees on $2n$ 
edges to Dyck $n$-paths: a worm crawls counterclockwise around the 
tree starting just left of the root and when an edge is traversed for 
the first time, records an upstep if the edge is left-leaning and a downstep 
if it is right-leaning. This bijection carries deep interior vertices to $DUU$s where deep interior means ``neither a 
leaf nor adjacent to a leaf''. Hence the recurrence also counts 
full binary trees on $2n$ edges by number 
of deep interior vertices. The interpretation for $m=2$ below is 
analogous to this one.

\vspace*{3mm}

{\Large \textbf{3 \quad Aligned cubic trees}  }

It is well known that there are $C_{n}$ (Catalan number) full binary 
trees on $2n$ edges. Considered as a graph, the root is the only vertex 
of degree 2 when $n \ge 1$. To remedy this, add a vertical planting 
edge to the root and transfer the root to the new vertex. Now every 
vertex has degree 1 or 3 and, throughout this paper, we will refer to a 
vertex of degree 3 as a \emph{node} and of degree 1 as a \emph{leaf}. 
Thus our planted tree of $2n+1$ edges has $n$ nodes. Leave the root 
edge pointing South and align the other edges so that all three 
angles at each node are $120^{\circ}$, lengthening edges as needed to 
avoid self intersections.  Rotate these objects through 
multiples of  $60^{\circ}$ to get all ``rooted aligned cubic'' 
trees on $n$ nodes ($6C_{n}$ of them, since the edge from the root is 
no longer restricted to point South but may point in any of 6 
directions). Now erase the root on each to get all (unrooted) aligned 
cubic trees on $n$ nodes ($\frac{6}{n+2}C_{n}$ of them since for 
each, a root could be placed on any of its $n+2$ leaves). Thus two 
drawings of an aligned cubic tree are equivalent if they 
differ only by translation and length of edges.
This is interpretation (iii) of the second order super-Catalan numbers 
mentioned in the Introduction. For short, we will refer to an aligned 
cubic tree on $n$ nodes simply as an $n$-ctree (c for cubic).

More concretely, a rooted $n$-ctree can be coded as a pair $(r,u)$ 
with $r$ an integer mod 6 and $u$ a nonnegative integer sequence of length 
$n+2$. The integer $r$ gives the angle (in multiples of $60^{\circ}$) 
from the direction South counterclockwise to the direction of the edge 
from the root. The sequence $u=(u_{i})_{i=1}^{n+2}$ gives the ``distance'' between 
successive leaves: traverse the tree in preorder (a worm crawls counterclockwise 
around the tree starting at a point just right of the root when 
looking from the root along the root edge). Then $u_{i}=v_{i}-2$ 
where $v_{i}\ (\ge 2)$ is the number of edges traversed between the 
$i$th leaf and the next one. 
For example, the sketched 3-ctree when 
rooted at $A$ is coded by $\left(2,(3,0,1,1,1,0)\right)$ and when 
rooted at $B$ is coded by $\left(1,(0,1,1,1,0,3)\right)$. 
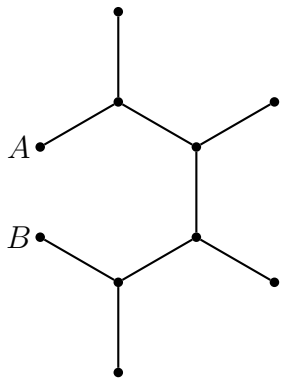
\begin{center}
\psset{xunit=.6cm,yunit=.6cm}
\pspicture*(-3,-1)(4,9)
     \rput{*0}(-2.2,3){$B$}
     \rput{*0}(-2.2,5){$A$}  
     \dotnode(0,8){a1}
     \dotnode(0,6){b1}
     \dotnode(3.46,6){b2}
     \dotnode(-1.73,5){c1}
     \dotnode(1.73,5){c2}
     \dotnode(-1.73,3){d1}
     \dotnode(1.73,3){d2}
     \dotnode(0,2){e1}
     \dotnode(3.46,2){e2}
     \dotnode(0,0){f1}
      \ncline{b1}{a1}
      \ncline{b1}{c1}
      \ncline{b1}{c2}
      \ncline{b2}{c2}
      \ncline{d2}{c2}
      \ncline{d2}{e1}
      \ncline{d2}{e2}
      \ncline{e1}{d1}
      \ncline{e1}{f1}
\endpspicture
\end{center}
In general, 
if a ctree rooted at a given leaf is coded by 
$\left(r,(u_{1},u_{2},\ldots,u_{n+1},u_{n+2})\right)$ then, when 
rooted at the next leaf in preorder, it is coded by $( 
(2+r-u_{1})$\,mod\,6$,(u_{2},u_{3},\ldots,u_{n+2},u_{1}))$. 
Repeating this $n+2$ times all told rotates $u$ back to itself and 
gives ``$r$''\,$=2n+4+r-\sum_{i=1}^{n+2}u_{i}$ mod 6. Since 
$\sum_{i=1}^{n+2}u_{i}$ is necessarily $=2n-2$, we are, as expected, 
back to the original coding sequence. 

An ordinary (planted) full binary tree is coded by $(0,u)$ and so there 
are $C_{n}$ coding sequences of length $n+2$. They can be generated as 
follows. A ctree can be built up by successively adding two edges to a 
leaf to turn it into a node. The effect this has on the coding sequence is 
to take two consecutive entries $u_{i},u_{i+1}$ (subscripts 
modulo $n+2$) and replace them by the three entries $u_{i}+1,\,0,\,u_{i+1}+1$. The 
1-ctree has coding sequence (0,0,0). The 2-ctree coding sequences are 
(1,0,1,0) and (0,1,0,1), and so on. Reversing this procedure gives a 
fast computational method to check if a given $u$ is a coding sequence 
or not. For example, successively pruning the first 0, $11210230 
\rightarrow 1120130 \rightarrow 111030 \rightarrow 11020 \rightarrow 1010 
\rightarrow 000$ is indeed a coding sequence.
However, we will work with the graphical depiction of a ctree.

The $n$-ctrees for $n=0,1,2$ are shown below. Note that since 
edges have a fixed non-horizontal direction, we can distinguish a top 
and bottom vertex for each edge. 
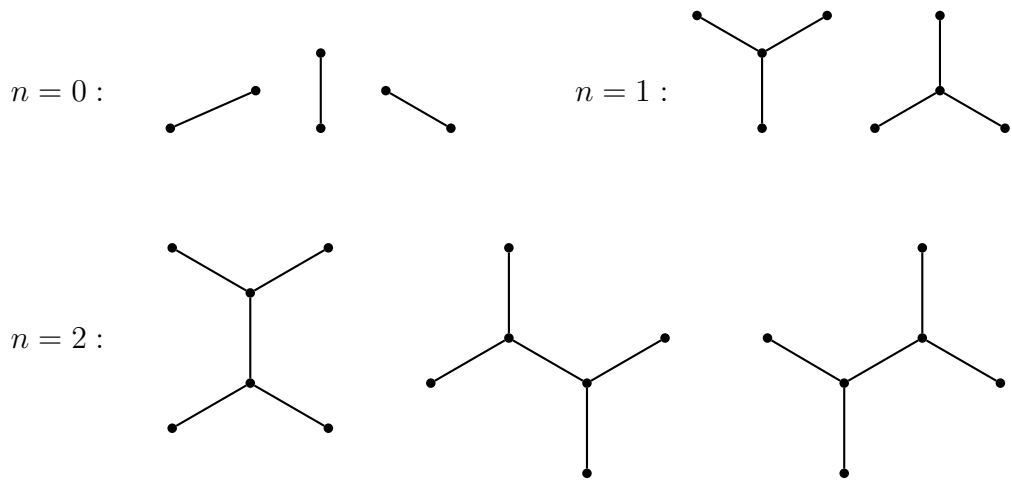
\begin{center}
\psset{xunit=.5cm,yunit=.5cm}
\pspicture*(-8,4)(8,9)
     \rput{*0}(-5,6){$n=0:$}    
     \dotnode(.27,6){A}
     \dotnode(2,7){B}
     \dotnode(3.73,6){C}
     \dotnode(-2,5){D}
     \dotnode(2,5){E}
     \dotnode(5.46,5){F}
      \ncline{A}{D}
      \ncline{B}{E}
      \ncline{C}{F}
\endpspicture
\pspicture*(-6,4)(8,9)
      \rput{*0}(-4,6){$n=1:$}
     \dotnode(-2,8){G}
     \dotnode(1.46,8){H}
     \dotnode(4.46,8){I}
     \dotnode(-.27,7){J}
     \dotnode(4.46,6){K}
     \dotnode(-.27,5){L}
     \dotnode(2.73,5){M}
     \dotnode(6.19,5){N}     
      \ncline{J}{G}
      \ncline{J}{H}
      \ncline{J}{L}
      \ncline{K}{I}
      \ncline{K}{M}
      \ncline{K}{N}
\endpspicture

\end{center}
\begin{center}
\psset{xunit=.6cm,yunit=.6cm}
\pspicture*(-12,-1)(12,5.5)
     \rput{*0}(-10,3){$n=2:$}    
     \dotnode(-7.46,5){a1}
     \dotnode(-4,5){a2}
     \dotnode(0,5){a3}
     \dotnode(9.16,5){a4}
     
     \dotnode(-5.73,4){b1}
     \dotnode(0,3){b2}
     \dotnode(3.46,3){b3}
     \dotnode(5.73,3){b4}
     \dotnode(9.16,3){b5}
     
     \dotnode(-5.73,2){c1}
     \dotnode(-1.73,2){c2}
     \dotnode(1.73,2){c3}
     \dotnode(7.43,2){c4}
     \dotnode(10.89,2){c5}
     
     \dotnode(-7.46,1){d1}
     \dotnode(-4,1){d2}
     \dotnode(1.73,0){d3}
     \dotnode(7.43,0){d4}
      \ncline{a1}{b1}
      \ncline{a2}{b1}
      \ncline{c1}{b1}
      \ncline{c1}{d1}
      \ncline{c1}{d2}
      
      \ncline{a3}{b2}
      \ncline{c2}{b2}
      \ncline{c3}{b2}
      \ncline{c3}{d3}
      \ncline{c3}{b3}
      
      \ncline{a4}{b5}
      \ncline{c4}{b5}
      \ncline{c5}{b5}
      \ncline{c4}{d4}
      \ncline{c4}{b4}

\endpspicture

\end{center}
It is convenient to introduce some further terminology. Recall a node 
is a vertex of degree 3. A node is \emph{hidden, exposed, naked} or 
\emph{stark naked} according as its 3 neighbors include $0,1,2$ or 3 
leaves. Thus a deep interior vertex is just a hidden node. A 0-ctree 
has no nodes. Only a 1-ctree has a stark naked node, and hidden nodes 
don't occur until $n\ge 4$. For $n\ge 2$, an $n$-ctree containing $k$ 
hidden nodes has $k+2$ naked nodes and hence $n-2k-2$ exposed nodes. The 
terms right and left can be ambiguous: we always use right and left 
relative to travel \emph{from} a specified vertex or edge. Thus vertex $B$ 
below is left (not right!) travelling from vertex $A$.

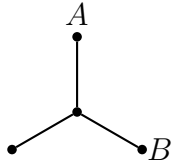
\vspace*{-30mm}

\hspace*{40mm}\pspicture*(-2,0)(2,6.5)
\psset{xunit=.5cm,yunit=.5cm}   
     \dotnode(0,5){a3}    
     \dotnode(0,3){b2}    
     \dotnode(-1.73,2){c2}
     \dotnode(1.73,2){c3}

      \ncline{a3}{b2}
      \ncline{c2}{b2}
      \ncline{c3}{b2}

      \rput(0,5.6){$A$}
      \rput(2.2,2){$B$}
      
\endpspicture \\
Each $n$-ctree has a unique center, either an edge or a node, defined 
as follows. For $n=0$, it is the (unique) edge in the ctree. For $n=1$, 
it is the (unique) node in the ctree. For $n\ge 2$, delete the leaves 
(and incident edges) adjacent to each naked node, thereby reducing the 
number of nodes by at least 2. Repeat until the $n=0$ or 1 definition 
applies. Equivalently, define the depth of a node in a ctree to be 
the length (number of edges) in the shortest path from the node to a 
naked node. Then there are either one or two nodes of maximal depth; 
if one, it is the center and if two, they are adjacent and the edge 
joining them is the center. 

\vspace*{3mm}
  
    {\Large \textbf{4 \quad (\emph{n\,,\,k})-Configurations}  }

There are $\binom{n-2}{2k}2^{n-2-2k}u_{k}$ configurations formed in 
the following way. Start with a $k$-ctree---$u_{k}$ choices. Break a 
strip of $n-2-2k$ squares---$\underbrace{\Box\!\Box\!\Box\!\Box \ldots 
\Box\!\Box}_{n-2-2k}$---into 
$2k+1$ (possibly empty) substrips, one for each of the $2k+1$ edges 
in the $k$-ctree---$\binom{(n-2-2k)+(2k+1)-1}{n-2-2k}=\binom{n-2}{2k}$ choices. 
Mark each square $L\,(=$ 
left) or $R\,(=$ right)---$2^{n-2-2k}$ choices.

Actually, this is not quite what we want. Perform one little 
tweaking: if there is a center edge and it has an 
\emph{odd-length} strip of squares, mark the first square $T\,(=$ 
top) or $B\,(=$ bottom) instead of $L$ or $R$. So a configuration 
might look as follows ($n=12,\,k=2$, empty strips not shown).

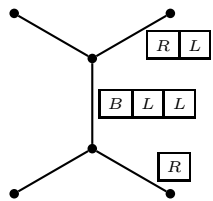
\begin{center}
\psset{xunit=.6cm,yunit=.6cm} 
\pspicture*(-2,-1)(3,5)
    
 \dotnode(-1.73,4){a}
      \dotnode(1.73,4){b}
      \dotnode(0,3){c}
      \dotnode(0,1){d}
       \dotnode(-1.73,0){e}
      \dotnode(1.73,0){f}
      \ncline{a}{c}
      \ncline{c}{b}
      \rput{*0}(2.1,3.3){{\tiny \fbox{$R$}\fbox{$L$}   } }
      \ncline{d}{e}
      \ncline{c}{d}
       \rput{*0}(1.4,2){{\tiny \fbox{$B$}\fbox{$L$}\fbox{$L$}   } }
      \ncline{d}{f}
       \rput{*0}(2.0,0.6){{\tiny \fbox{$R$}   } }
\endpspicture

\end{center}
{\Large \textbf{5 \quad Bijection}  }

Here is a bijection from $(n,k)$-configurations to $n$-ctrees with $k$ 
hidden nodes. The bijection produces the correspondences in the 
following table.
\begin{center}
\begin{tabular}{cc}
   \textbf{(\emph{n,k})-configuration} &   \textbf{\emph{n}-ctree}    \\
   leaf & naked node  \\
   node & hidden node  \\
    square & exposed node  \\
\end{tabular}
\end{center}

Roughly speaking, work outward from the center, turning a strip of 
$j$ labeled squares on an edge $AB$ into $j$ exposed nodes lying 
between $A$ and $B$.

First, for the center edge (if there is one), the procedure depends 
on whether it has an even or odd number of squares.

\textbf{Case Even}\quad Here, the center edge becomes an edge joining two exposed nodes as 
shown: the labels again 
indicate the 
$L/R$ status (travelling from the center edge) of the leaves associated 
with the exposed nodes.
The labels are applied from the bottom vertex subtree
$(H_{2})$ to the top one $(H_{1})$.
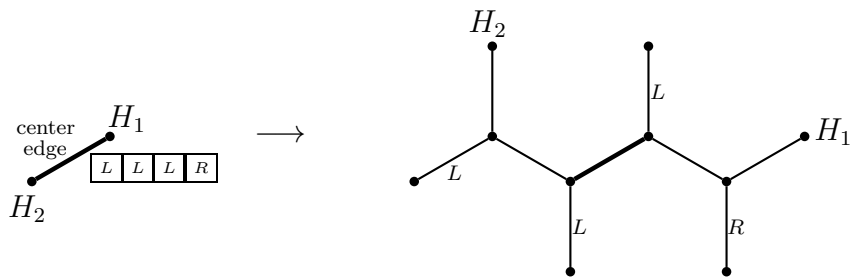
\begin{center}

\pspicture*(-4,-1)(4,4.5)
\psset{xunit=.6cm,yunit=.6cm}     
      \dotnode(-1,2){a}
      \dotnode(.73,3){b}

      \ncline[linewidth=.06]{a}{b}
      
      \rput{*0}(-0.6,3.2){{\scriptsize center }}
        \rput{*0}(-0.7,2.7){{\scriptsize edge}}
      
       \rput{*0}(-1.1,1.4){$H_{2}$}
       \rput{*0}(1.1,3.4){$H_{1}$}
       \rput{*0}(4.5,3){$\longrightarrow$}

       \rput{*0}(1.9,2.3){{\tiny \fbox{$L$}\fbox{$L$}\fbox{$L$}\fbox{$R$}   } }
     
\endpspicture
\pspicture*(-3,-1)(6,4.5)
\psset{xunit=.6cm,yunit=.6cm}     
     \rput{*0}(-2.56,5.5){$H_{2}$}
     \rput{*0}(5.1,3.1){$H_{1}$}
     
     \rput{*0}(-3.3,2.2){{\scriptsize $L$}}
     \rput{*0}(-.53,1){{\scriptsize $L$}}
     \rput{*0}(1.2,4){{\scriptsize $L$}}
      \rput{*0}(2.93,1){{\scriptsize $R$}}

     \dotnode(-2.46,5){a1}
     \dotnode(1,5){a2}

     \dotnode(-2.46,3){b1}
     \dotnode(1,3){b2}
     \dotnode(4.46,3){b3}

     \dotnode(-4.19,2){c1}
     \dotnode(-.73,2){c2}
     \dotnode(2.73,2){c3}

     \dotnode(-.73,0){d1}
     \dotnode(2.73,0){d2}
     
      \ncline{a1}{b1}
      \ncline{c2}{b1}
      \ncline{c1}{b1}

      \ncline{a2}{b2}
      \ncline{c2}{d1}
      \red{ \ncline[linewidth=.06]{c2}{b2}}
      \ncline{c3}{b2}
      \ncline{c3}{d2}
      \ncline{c3}{b3}
            
\endpspicture

\end{center}
\textbf{Case Odd}\quad Here, the center edge becomes a leaf edge. 
The first square indicates whether the top or bottom vertex becomes 
a leaf. Construct equal numbers of exposed nodes on each 
side of the non-leaf (here, top) vertex, using the $L/R$ designations 
to determine the leaves ($L/R$ relative to travel from the leaf and running, 
say, from the left branch to the right).
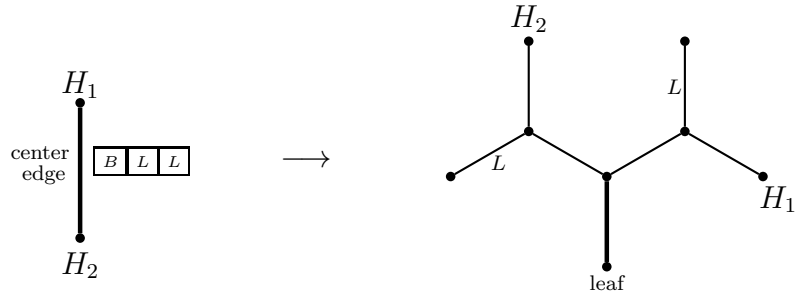
\begin{center}
    
\pspicture*(-5,0)(4,4)
\psset{xunit=.6cm,yunit=.6cm}         
      
      \dotnode(0,5.3){c}
      \dotnode(0,2.3){d}
     
       \rput{*0}(0.0,5.7){$H_{1}$}
       \rput{*0}(0.0,1.7){$H_{2}$}
       \rput{*0}(5,4){$\longrightarrow$}
       
       \rput{*0}(-0.8,4.2){{\scriptsize center }}
        \rput{*0}(-0.8,3.7){{\scriptsize edge}}
      
      \ncline[linewidth=.06]{c}{d}

       \rput{*0}(1.55,4){{\tiny \fbox{$B$}\fbox{$L$}\fbox{$L$}   } }

\endpspicture    
\pspicture*(-3,-1)(4,3.5)
\psset{xunit=.6cm,yunit=.6cm}         
      
      \dotnode(-1.73,5){a}
      \dotnode(1.73,5){b}
      \dotnode(-1.73,3){c}
      \dotnode(1.73,3){d}
      \dotnode(-3.46,2){e}
      \dotnode(0,2){f}
      \dotnode(3.46,2){g}
      \dotnode(0,0){h}
      
      \ncline{c}{a}
      \ncline{c}{e}
      \ncline{c}{f}
      \ncline{b}{d}
      \ncline{f}{d}
      \ncline{g}{d}
      \ncline[linewidth=.06]{f}{h}
      
       \rput{*0}(0,-.35){{\scriptsize leaf}}
       \rput{*0}(-2.4,2.3){{\scriptsize $L$}}
       \rput{*0}(1.5,4){{\scriptsize $L$}}
       \rput{*0}(-1.73,5.5){$H_{2}$}
       \rput{*0}(3.8,1.5){$H_{1}$}
     
\endpspicture

\end{center}
Decide the placement of the two subtrees $H_{1},H_{2}$, say the $H$ originally sitting 
at the vertex which is now a leaf goes at the end of the left branch 
from the leaf.

Next, for a non-center edge with $i\ge 0$ squares, identify its 
endpoint closest to the center, let $H_{0},H_{1},H_{2}$ be the 
subtrees as illustrated ($H_{0}$ containing the center), and insert $i$ exposed nodes 
as shown. The labels $L,L,R$ apply in order from 
$H_{1}$-$H_{2}$ to $H_{0}$ and indicate the 
$L/R$ status (travelling from the center) of the leaves associated 
with the exposed nodes.
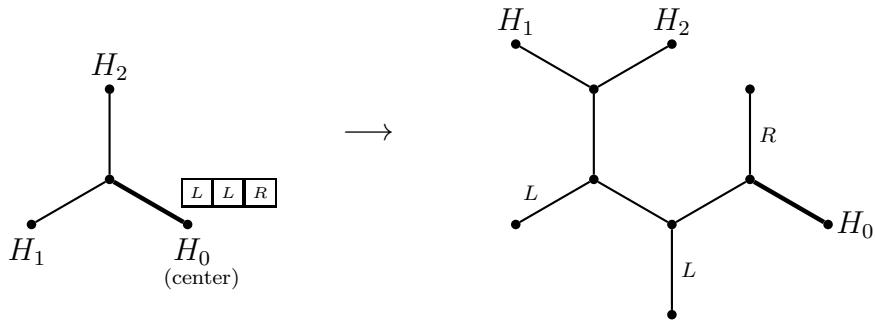
\begin{center}
\psset{xunit=.6cm,yunit=.6cm} 
\pspicture*(-14,-1)(7,7)
    
      \dotnode(-1.73,6){a1}
      \dotnode(1.73,6){a2}
      \dotnode(-10.73,5){b1}
      \dotnode(0,5){b2}
      \dotnode(3.46,5){b3}
      \dotnode(-10.73,3){c1}
      \dotnode(0,3){c2}
      \dotnode(3.46,3){c3}
      \dotnode(-12.46,2){d1}
      \dotnode(-9,2){d2}
      \dotnode(-1.73,2){d3}
      \dotnode(1.73,2){d4}
      \dotnode(5.19,2){d5}
      \dotnode(1.73,0){e1}

      \ncline{c1}{b1}
      \ncline{c1}{d1}
      \ncline[linewidth=.06]{c1}{d2}
      
      \ncline{a1}{b2}
      \ncline{a2}{b2}
      \ncline{c2}{b2}
      \ncline{c2}{d3}
      \ncline{c2}{d4}
      \ncline{c3}{d4}
      \ncline{e1}{d4}
      \ncline{c3}{b3}
      \ncline[linewidth=.06]{c3}{d5}

       \rput{*0}(-10.73,5.5){$H_{2}$}
       \rput{*0}(-12.56,1.4){$H_{1}$}
       \rput{*0}(-8.9,1.4){$H_{0}$}
        \rput{*0}(-8.6,0.8){{\scriptsize(center) }}
       \rput{*0}(-5,4){$\longrightarrow$}
       
       \rput{*0}(-1.73,6.5){$H_{1}$}
       \rput{*0}(1.73,6.5){$H_{2}$}
       \rput{*0}(5.8,2){$H_{0}$}
       
       \rput{*0}(2.1,1){{\scriptsize $L$}}
       \rput{*0}(3.85,4){{\scriptsize $R$}}
       \rput{*0}(-1.4,2.7){{\scriptsize $L$}}      
    
       \rput{*0}(-7.9,2.7){{\tiny \fbox{$L$}\fbox{$L$}\fbox{$R$}   } }

\endpspicture

\end{center}

Finally, turn the original leaves into naked nodes by adding two edges apiece.

We leave to the reader to verify that the resulting ctree has $n$ 
nodes of which $k$ are hidden, and that the original configuration 
can be uniquely recovered from this $n$-ctree by reversing the above 
procedure, working in from the leaves.

\end{document}